\documentclass[10pt,letterpaper]{article}
\usepackage{blindtext,titlefoot}
\usepackage{authblk}
\usepackage{marvosym}
\usepackage[utf8]{inputenc}
\usepackage[english]{babel}
\usepackage{amsmath}
\usepackage{amsfonts}
\usepackage{amssymb}
\usepackage{makeidx}
\usepackage{graphicx}
\usepackage{lmodern}
\usepackage{kpfonts}
\usepackage{dsfont}
\usepackage{cancel}
\usepackage{amsthm} 
\usepackage[left=2cm,right=2cm,top=2cm,bottom=2cm]{geometry}
\usepackage{subcaption}
\usepackage{multirow}
\usepackage{float}
\usepackage{url}
\usepackage{breakurl}
\usepackage[breaklinks]{hyperref}
\usepackage{multicol}
\usepackage{mathtools, cuted}
\usepackage{lipsum}
\usepackage{booktabs}
\usepackage{comment}
\usepackage{cite}

\usepackage{longtable}

\DeclareMathOperator{\re}{Re}
\DeclareMathOperator{\im}{Im}
\DeclareMathOperator{\sgn}{sgn}
\DeclareMathOperator{\tr}{tr}
\DeclareMathOperator{\card}{card}
\DeclareMathOperator{\Op}{O}
\DeclareMathOperator{\Opw}{W}
\DeclareMathOperator{\Ops}{S}
\DeclareMathOperator{\Opt}{T}
\DeclareMathOperator{\Opr}{R}
\DeclareMathOperator{\Conv}{Conv}
\DeclareMathOperator{\rnd}{Rnd}
\providecommand{\nset}[1]{
\mathbb{#1}
}
\providecommand{\set}[1]{
\left\{#1\right\}
}

\providecommand{\ifr}[5]{
{}^{#1}_{#2}{#3}_{#4}^{#5}
}
\providecommand{\gam}[1]{
\Gamma\left(#1 \right)
}

\providecommand{\norm}[1]{
\left\lVert #1 \right\rVert
}
\providecommand{\abs}[1]{
\left\lvert #1 \right\rvert
}
\providecommand{\ds}[1]{
\displaystyle #1
}
\providecommand{\der}[3]{
\dfrac{#1^{#3} }{ #1 #2^{#3}}
}

\providecommand{\rnd}[2]{
\hbox{Rnd}_#2\left(#1\right)
}
\newcommand{\mref}[2]{\textbf{#1 #2}}
\newtheorem{theorem}{ Theorem}[section]
\newtheorem{definition}[theorem]{Definition}
\newtheorem{proposition}[theorem]{Proposition}
\newtheorem{corollary}[theorem]{Corollary}
\newtheorem{example}{Example}

\usepackage{enumitem} 
\setlist[itemize]{noitemsep} 

\usepackage{titlesec} 
\titleformat{\section}[block]{\large\bfseries\scshape\centering}{\thesection.}{1em}{} 
\titleformat{\subsection}[block]{\large\bfseries\scshape\centering}{\thesubsection.}{1em}{}
\titleformat{\subsubsection}[block]{\large\bfseries\scshape\centering}{\thesubsubsection.}{1em}{} 

\title{\Huge\bfseries Sets of fractional operators and numerical estimation of the order of convergence of a family of fractional fixed point methods}

\author[,a]{A. Torres-Hernandez  \footnote{Email: anthony.torres@ciencias.unam.mx; Corresponding author; ORCID: 0000-0001-6496-9505}}
\affil[a]{Department of Physics, Faculty of Science - UNAM, Mexico}

\author[,b]{F. Brambila-Paz \footnote{Email: fernandobrambila@gmail.com; ORCID: 0000-0001-7896-6460}}
\affil[b]{Department of Mathematics, Faculty of Science - UNAM, Mexico}

\date{}

\begin{document}

\maketitle


\begin{abstract}

Considering the large number of fractional operators that exist, and since it does not seem that their number will stop increasing soon at the time of writing this paper, it is presented for the first time, as far as the authors know, a simplified and compact way to work the fractional calculus through the classification of fractional operators using sets.  This new way of working with fractional operators, which may be called as fractional calculus of sets, allows to generalize objects of the conventional calculus such as tensor operators, the diffusion equation, the heat equation, the Taylor series of a vector-valued function, and the fixed point method in several variables which allows to generate the method known as the fractional fixed point method.  It is also shown that each fractional fixed point method that generates a convergent sequence has the ability to generate an uncountable family of fractional fixed point methods that generate convergent sequences. So, it is shown one way to estimate numerically the mean order of convergence of any fractional fixed point method in a region $\Omega$ through the problem of determining the critical points of a scalar function, and it is shown how to construct a hybrid fractional iterative method to determine the critical points of a scalar function.

\end{abstract}

\section{Introduction}

A fractional derivative is an operator that generalizes the ordinary derivative, in the sense that if

\begin{eqnarray*}
\dfrac{d^\alpha}{d x^\alpha},
\end{eqnarray*}

denotes the differential of order $ \alpha $, then $\alpha$ may be considered a parameter, with $\alpha\in \nset{R}$, such that the first derivative corresponds to the particular case $\alpha=1$. On the other hand, a fractional differential equation is an equation that involves at least one differential operator of order $ \alpha $, with $(n-1)< \alpha \leq n$ for some positive integer $n$, and it is said to be a differential equation of order $\alpha$ if this operator is the highest order in the equation. Analogously, a fractional partial differential equation is an equation that involves at least one differential operator of order $ \alpha $, which in general are usually partial derivatives of order $ \alpha $, that is,

\begin{eqnarray*}
\dfrac{\partial^\alpha}{\partial t^\alpha}, & \dfrac{\partial^\alpha}{\partial x^\alpha}, &\dfrac{\partial^\alpha}{\partial y^\alpha}.
\end{eqnarray*}

The fractional operators have many representations, but one of their fundamental properties is that they allow retrieving the results of conventional calculus when $\alpha \to n$. So, considering a scalar function $h: \nset{R}^m \to \nset{R}$ and the canonical basis of $\nset{R}^m$ denoted by $\set{\hat{e}_k}_{k\geq 1}$, and it is possible to define the following fractional operators of order $\alpha$ using Einstein notation

\begin{eqnarray}
o_x^\alpha h(x):=\hat{e}_k o_k^\alpha h(x),
\end{eqnarray}

then denoting by  $\partial_k^n$ the partial derivative of order $n$ applied with respect to the $k$-th component of the vector $x$, using the previous operator it is possible to define the following set of fractional operators

\begin{eqnarray}\label{eq:0}
\Op_{x,\alpha}^n(h):=\set{o_x^\alpha \ : \ \exists \partial_k^n h(x)  \  \mbox{ and } \ \lim_{\alpha \to n}o_k^\alpha h(x)=\partial_k^n h(x) \ \forall k\geq 1 },
\end{eqnarray}

which may be proved to be a nonempty set through the following set of fractional operators

\begin{eqnarray}
\set{
o_x^\alpha \ : \ o_k^\alpha h(x)= \left(\dfrac{\alpha}{n} \partial_k^n  + \left(1-\dfrac{\alpha}{n} \right) \partial_k^\alpha \right)h(x) \ \mbox{ and } \   \lim_{ \alpha \to n}\partial_k^\alpha h(x) \neq \partial_k^n h(x) \ \forall k\geq 1},
\end{eqnarray}

and which may be considered as a generating set of sets of \textbf{fractional tensor operators}. For example, considering $\alpha,n\in \nset{R}^d$ with $\alpha=\hat{e}_k[\alpha]_k$ and $n=\hat{e}_k [n]_k$, it is possible to define the following set of fractional tensor operators

\begin{eqnarray}\label{eq:00}
\Op_{x,\alpha}^{n}(h):=\set{o_x^\alpha \ : \ o_x^\alpha \in \Op_{x,[\alpha]_1}^{[n]_1}(h) \times \Op_{x,[\alpha]_2}^{[n]_2}(h)\times \cdots \times \Op_{x,[\alpha]_d}^{[n]_d}(h) }, 
\end{eqnarray}

therefore, considering a function $h: \nset{R}^m \times \nset{R}_{\geq 0} \to \nset{R}$, as well as the vectors $\alpha,n\in \nset{R}^3$ with $\alpha=\hat{e}_k[\alpha]_k$ and $n=\hat{e}_k[n]_k $,  it is possible to combine the sets \eqref{eq:0} and \eqref{eq:00} to define new sets of fractional operators related to the theory of differential equations, as shown with the following set

\begin{eqnarray}
\Opw_{t,x,\alpha}^{n}(h):=\set{w_{t,x}^\alpha=o_t^{[\alpha]_1}-\tr\left(o_{x}^{([\alpha]_2,[\alpha]_3)}\right)\ : \ o_t^{[\alpha]_1} \in \Op_{t,[\alpha]_1}^{[n]_1}(h) \ \mbox{ and } \  o_{x}^{([\alpha]_2,[\alpha]_3)}\in \Op_{x,([\alpha]_2,[\alpha]_3)}^{([n]_2,[n]_3)}(h) \  },
\end{eqnarray}

where $\tr(\cdot)$ denotes the trace of a matrix. So, denoting the Laplacian operator by $\nabla^2$, it is possible to obtain the following results:

\begin{eqnarray}
\mbox{If } w_{t,x}^\alpha\in \Opw_{t,x,\alpha}^{n}(h) \mbox{ with }n=(1,1,1) \ \Rightarrow \ \lim_{\alpha \to n}w_{t,x}^\alpha h(x,t) = \left(\partial_t- \nabla^2 \right)h(x,t),\vspace{0.1cm}\\
\mbox{If } w_{t,x}^\alpha\in \Opw_{t,x,\alpha}^{n}(h) \mbox{ with }n=(2,1,1) \ \Rightarrow \ \lim_{\alpha \to n}w_{t,x}^\alpha h(x,t) = \left(\partial_t^2- \nabla^2 \right)h(x,t),
\end{eqnarray}

which may generalize the diffusion equation and the heat equation respectively. To finish this section, it is necessary to mention that the applications of fractional operators have spread to different fields of science such as finance \cite{safdari2015radial,torres2021blackscholes}, economics \cite{traore2020model,torres2020nonlinear}, number theory through the Riemann zeta function \cite{guariglia2021fractional,torres2021zeta} and in engineering with the study for the manufacture of hybrid solar receivers \cite{de2021fractional,torres2020reduction,torres2020fracpseunew,torres2021codefracpseudo}. It should be mentioned that there is also a growing interest in fractional operators and their properties for the solution of nonlinear systems \cite{wang2021derivative,gdawiec2019visual,akgul2019fractional,torres2021fracnewrap,torres2020fracsome,torres2021fracnewrapaitken}, which is a classic problem in mathematics, physics and engineering, which consists of finding the set of zeros of a function $f:\Omega \subset \nset{R}^n \to \nset{R}^n$, that is,

\begin{eqnarray}\label{eq:1-001}
\set{\xi \in \Omega \ : \ \norm{f(\xi)}=0},
\end{eqnarray}

where $\norm{ \ \cdot \ }: \nset{R}^n \to \nset{R}$ denotes any vector norm, or equivalently

\begin{eqnarray}
\set{\xi \in \Omega \ : \ [f]_k(\xi)=0 \ \forall k\geq 1}.
\end{eqnarray}

where $[f]_k: \nset{R}^n \to \nset{R}$ denotes the $k$-th component of the function $f$. This paper presents a simplified and compact way to work the fractional calculus through the classification of fractional operators using sets.  This new way of working with fractional operators allows to generalize objects of the conventional calculus such as tensor operators, the diffusion equation, the heat equation, the Taylor series of a vector-valued function, and the fixed point method in several variables which allows to generate the method known as the fractional fixed point method. It is also shown that each fractional fixed point method that generates a convergent sequence has the ability to generate an uncountable family of fractional fixed point methods that generate convergent sequences. It is shown one way to estimate numerically the mean order of convergence of any fractional fixed point method in a region $\Omega$ through the problem of determining the critical points of a scalar function, and it is shown how to construct a hybrid fractional iterative method to determine the critical points of a scalar function.

\section{Fixed Point Method}

Let  $\Phi:\nset{R}^n \to \nset{R}^n$ be a function. It is possible to build a sequence $\set{x_i}_{i\geq 1}$  by defining the following iterative method

\begin{eqnarray}\label{eq:2-001}
x_{i+1}:=\Phi(x_i), & i=0,1,2,\cdots,
\end{eqnarray}

if it is fulfilled that $x_i\to \xi\in \nset{R}^n$ and if the function $\Phi$ is continuous around $\xi$, we obtain that

\begin{eqnarray}\label{eq:2-002}
\ds \xi=\lim_{i\to \infty}x_{i+1}=\lim_{i\to \infty}\Phi(x_i)=\Phi\left(\lim_{i\to \infty}x_i \right)=\Phi(\xi),
\end{eqnarray}

the above result is the reason by which the method \eqref{eq:2-001} is known as the \textbf{fixed point method}. Furthermore, the function $\Phi$ is called an \textbf{iteration function}.  Before continuing, it is necessary to define the order of convergence of an iteration function $\Phi$ \cite{plato2003concise,torres2021fracnewrapaitken}:

\begin{definition}
Let $ \Phi: \Omega \subset \nset{R}^ n \to \nset{R}^ n $ be an iteration function with a fixed point $ \xi \in \Omega $. So, the method \eqref{eq:2-001} is called  (locally) convergent, with an \textbf{order of convergence} of  order (at least) $p $ (with $p\geq 1 $), if there exist $ \delta> 0 $  and a non-negative constant $ C $ (with $ C <1 $ if $ p = 1 $), such that for any initial value $ x_0 \in B (\xi; \delta) $ it is fulfilled that

\begin{eqnarray}\label{eq:2-0021}
\norm{x_{i+1}-\xi}\leq C \norm{x_i-\xi}^p, & i=0,1,2,\cdots,
\end{eqnarray}

where $ C $ is called \textbf{convergence factor}.

\end{definition}

The order of convergence is usually related to the speed at which the sequence generated by an iteration function $ \Phi $ converges. It is necessary to mention that for the case $ p = 1 $ it is said that the function $ \Phi $ has an \textbf{order of convergence (at least) linear}, and for the case $p=2$  it is said that the function $ \Phi $ has an \textbf{order of convergence (at least) quadratic}. From the previous definition the following proposition is obtained:

\begin{proposition}\label{prop:2-001}
Let $\Phi:\nset{R}^n \to \nset{R}^n$ be an iteration function that defines a sequence $\set{x_i}_{i\geq 1}$ such that $x_i\to \xi\in \nset{R}^n$. So, if $\Phi$ has an order of convergence of order (at least) $p$  in $B(\xi;\delta)$, there exists a non-negative constant $K=K(C)$, such that for all values of the sequence $\set{x_i}_{i\geq 1}$ it is fulfilled that

\begin{eqnarray}
\norm{x_{i+1}-x_{i}}\leq K \norm{x_{i}-x_{i-1}}^p, & i=0,1,2,\cdots,
\end{eqnarray}

where $\norm{x_{-1}}:=0$.

\begin{proof}
Considering that $\Phi$ defines a sequence $\set{x_i}_{i\geq 1}$ and that it has an order of convergence of order (at least) $p$, it is possible to obtain the following inequality

\begin{eqnarray*}
\norm{x_{i+1}-x_{i}}\leq C \big{(} \norm{x_{i}-\xi}^p+\left(\norm{x_{i}-\xi}+\norm{x_{i}-x_{i-1}}\right)^p \big{)},
\end{eqnarray*}

as a consequence

\begin{eqnarray*}
\norm{x_{i+1}-x_{i}}\leq 2C \left(\norm{x_{i}-\xi}+\norm{x_{i}-x_{i-1}}\right)^p ,
\end{eqnarray*}

and since $x_i\to \xi$, there exists a positive constant $c$ such that

\begin{eqnarray*}
\norm{x_{i+1}-x_{i}} \leq 2cC \norm{x_{i}-x_{i-1}}^p = K \norm{x_{i}-x_{i-1}}^p.
\end{eqnarray*}

\end{proof}

\end{proposition}

From the previous proposition the following theorem is obtained:

\begin{theorem}\label{teo:2-001}
Let $\Phi:\nset{R}^n \to \nset{R}^n$ be an iteration function that defines a sequence $\set{x_i}_{i\geq 1}$ such that $x_i\to \xi\in \nset{R}^n$. So, if $\Phi$ has an order of convergence of order (at least) $p$  in $B(\xi;\delta)$, there exists a value $m\in \nset{N}$ such that for all subsuccession $\set{x_i}_{i\geq m}\in  B(\xi;1/2)$ that fulfills the following condition

\begin{eqnarray*}
\norm{x_{i+2}-x_{i+1}}\leq K \norm{x_{i+1}-x_{i}}^p, &  \forall i\geq m,
\end{eqnarray*}

there exist $\delta_K=\delta_K(C)>0$ and a sequence of values $P_i$ given by the following expression

\begin{eqnarray}\label{eq:2-0022}
P_i:=\dfrac{\log\left(\norm{x_i-x_{i-1}} \right)}{\log\left(\norm{x_{i-1}-x_{i-2}} \right)}, 
\end{eqnarray}

such that $\set{P_i}_{i\geq m+2}\in B(p;\delta_K)$.

\begin{proof}

Considering that $\Phi$ defines a sequence $\set{x_i}_{i\geq 1}$ and that it has an order of convergence of order (at least) $p$, from the \mref{Proposition}{\ref{prop:2-001}} it is possible to obtain the following inequality

\begin{eqnarray*}
\log\left(\norm{x_{i+2}-x_{i+1}}\right)-p \log\left( \norm{x_{i+1}-x_{i}} \right)\leq \log(K),
\end{eqnarray*}

assuming that there exists a subsequence $\set{x_i}_{i\geq m}\in B(\xi;1/2)$, then $\log\left(\norm{x_{i+1}-x_{i}}\right)<0 \ \forall i\geq m$. So, if the subsequence $\set{x_i}_{i\geq m}$ fulfills the above inequality

\begin{eqnarray*}
\dfrac{\log(K)}{\log\left( \norm{x_{i+1}-x_{i}} \right)}\leq \dfrac{\log\left(\norm{x_{i+2}-x_{i+1}}\right)}{\log\left( \norm{x_{i+1}-x_{i}} \right)}-p ,
\end{eqnarray*}

considering that $x\leq \abs{x} \ \forall x\in \nset{R}$, there exists a positive constant $c$ such that

\begin{eqnarray*}
\dfrac{\log(K)}{\log\left( \norm{x_{i+1}-x_{i}} \right)}\leq \abs{ \dfrac{\log\left(\norm{x_{i+2}-x_{i+1}}\right)}{\log\left( \norm{x_{i+1}-x_{i}} \right)}-p} \leq c\abs{ \dfrac{\log(K)}{\log\left( \norm{x_{i+1}-x_{i}} \right)}},
\end{eqnarray*}

and since $K=K(C)$, there exists  a positive value $\delta_K=\delta_K(C)$ such that the sequence $\set{P_i}_{i\geq m+2}\in B(p;\delta_K)$.

\end{proof}

\end{theorem}

From the previous theorem the following corollary is obtained:

\begin{corollary}
Let $\Phi:\nset{R}^n \to \nset{R}^n$ be an iteration function that defines a sequence $\set{x_i}_{i\geq 1}$ such that $x_i\to \xi\in \nset{R}^n$. So, if $\Phi$ has an order of convergence of order (at least) $p$  in $B(\xi;1/2)$, for some $m\in \nset{N}$ there exists a sequence $\set{P_i}_{i\geq m}\in B(p;\delta_K)$ that fulfills the following condition

\begin{eqnarray*}
\lim_{i\to \infty}P_i \to p,
\end{eqnarray*}

and therefore, there exists at least one value $k\geq m$ such that

\begin{eqnarray}
\abs{\and{P_k-p}}\leq \epsilon.
\end{eqnarray}

\end{corollary}

The previous corollary allows estimating numerically the order of convergence of an iteration function $\Phi$ that generates at least one convergent sequence $\set{x_i}_{i\geq 1}$. On the other hand, the following corollary allows characterizing the order of convergence of an iteration function $ \Phi $ through its \textbf{Jacobian matrix} $\Phi^{(1)}$ \cite{torres2021fracnewrapaitken}:

\begin{corollary}\label{cor:2-001}
Let $\Phi:\nset{R}^n \to \nset{R}^n$ be an iteration function. If $\Phi$ defines a sequence $\set{x_i}_{i\geq 1}$ such that $x_i\to \xi\in \nset{R}^n$. So, $\Phi$ has an order of convergence of order (at least) $p$  in $B(\xi;\delta)$, where

\begin{eqnarray}\label{eq:2-003}
p:=\left\{
\begin{array}{cc}
1 , &\ds \mbox{if } \lim_{x\to \xi}\norm{\Phi^{(1)}(x)}\neq 0  \vspace{0.2cm}\\
2, &\ds \mbox{if } \lim_{x\to \xi}\norm{\Phi^{(1)}(x)}= 0  \\
\end{array}\right. .
\end{eqnarray}

\end{corollary}

\section{Riemann-Liouville Fractional Operators}

One of the fundamental operators of fractional calculus is the operator \textbf{Riemann-Liouville fractional integral},  which is defined as follows \cite{hilfer00,oldham74}

\begin{eqnarray}
\ifr{}{a}{I}{x}{\alpha}f(x):=\dfrac{1}{\gam{\alpha}}\int_a^x (x-t)^{\alpha-1}f(t)dt,
\end{eqnarray}

which is a fundamental piece to build the operator \textbf{Riemann-Liouville fractional derivative},  which is defined as follows \cite{hilfer00,kilbas2006theory}

\begin{eqnarray}\label{eq:3-001}
\ifr{}{a}{D}{x}{\alpha}f(x) := \left\{
\begin{array}{cc}
\ds \ifr{}{a}{I}{x}{-\alpha}f(x), &\mbox{if }\alpha<0 \vspace{0.1cm}\\  
\ds \dfrac{d^n}{dx^n}\left( \ifr{}{a}{I}{x}{n-\alpha}f(x)\right), & \mbox{if }\alpha\geq 0
\end{array}
\right., 
\end{eqnarray}

where  $ n = \lceil \alpha \rceil$ and $\ifr{}{a}{I}{x}{0}f(x):=f(x)$. Applying the operator \eqref{eq:3-001} with $a=0$ to the  function $ x^{\mu} $, with $\mu> -1$, we obtain the following result \cite{torres2021fracnewrapaitken}:

\begin{eqnarray}\label{eq:3-002}
\ifr{}{0}{D}{x}{\alpha}x^\mu = 
 \dfrac{\gam{\mu+1}}{\gam{\mu-\alpha+1}}x^{\mu-\alpha}, & \alpha\in \nset{R}\setminus \nset{Z}.
\end{eqnarray}

\section{Fractional Fixed Point Method}

Let $\nset{N}_0$ be the set $\nset{N}\cup\set{0}$, if $\gamma \in \nset{N}_0^m$ and $x\in \nset{R}^m$, then it is possible to define the following multi-index notation

\begin{eqnarray}
\left\{
\begin{array}{c}
\begin{array}{ccc}
 \gamma!:= \ds\prod_{k=1}^m [\gamma]_k !,& \abs{\gamma}:= \ds \sum_{k=1}^m [\gamma]_k\vspace{0.1cm}, &  x^\gamma:= \ds \prod_{k=1}^m [x]_k^{[\gamma]_k}
\end{array} \\
\der{\partial}{x}{\gamma}:= \dfrac{\partial^{\abs{\gamma}}}{\partial [x]_1^{[\gamma]_1}\partial [x]_2^{[\gamma]_2}\cdots \partial [x]_m^{[\gamma]_m} }
\end{array}\right. .
\end{eqnarray}

So, considering a function $h: \Omega \subset \nset{R}^m \to \nset{R}$, it is possible to define the following set of fractional operators

\begin{eqnarray}
\Ops_{x,\alpha}^{n,\gamma}(h):=\set{s_x^{\alpha \gamma}=s_x^{\alpha\gamma}\left( o_x^\alpha \right) \ : \ o_x^\alpha \in \Op_{x,\alpha}^{k\abs{\gamma}}(h) \ \mbox{ and } \ s_x^{\alpha \gamma}h(x):=o_1^{\alpha[\gamma]_1}o_2^{\alpha [\gamma]_2}\cdots o_m^{\alpha [\gamma]_m}h(x) \ \forall k,\abs{\gamma}\leq n },
\end{eqnarray}

from which it is possible to obtain the following result:

\begin{eqnarray}
\mbox{If }s_x^{\alpha \gamma}\in \Ops_{x,\alpha}^{1,\gamma}(h) \ \Rightarrow \ \lim_{\alpha \to 1}s_x^{\alpha \gamma}h(x)=\der{\partial}{x}{\gamma}h(x) \ \forall \abs{\gamma}\leq 1,
\end{eqnarray}

and as a consequence, considering a function $h: \Omega \subset \nset{R}^m \to \nset{R}^m$, it is possible to define the following set of fractional operators

\begin{eqnarray}
{}_m\Ops_{x,\alpha}^{n,\gamma}(h):=\set{ s_x^{\alpha \gamma} \ : \ s_x^{\alpha \gamma} \in \Ops_{x,\alpha}^{n,\gamma}\left( [h]_k \right) \ \forall k \leq m}.
\end{eqnarray}

On the other hand, using little-o notation it is possible to obtain the following result:

\begin{eqnarray}
\mbox{If }x\in B(a;\delta)  \ \Rightarrow \ \lim_{ x\to a  } \dfrac{o\left((x-a)^\gamma \right)}{(x-a)^\gamma}\to 0 \ \forall \abs{\gamma}\geq 1,
\end{eqnarray}

with which it is possible to define the following set of functions

\begin{eqnarray}
\Opr_{\alpha\gamma}^{n}(a):=\set{r_{\alpha\gamma}^n  \ :  \ \dfrac{\norm{r_{\alpha\gamma}^n(x-a)}}{\norm{x-a}^n}\leq \dfrac{o\left(\norm{x-a}^n\right)}{\norm{x-a}^n} \ \forall \abs{\gamma}\geq n \ \mbox{ and } \ \forall x\in B(a;\delta)}.
\end{eqnarray}

So, considering the previous set and some $B(a;\delta)\ \subset \Omega$, it is possible to define the following sets of fractional operators

\begin{eqnarray}
{}_m\Opt_{x,\alpha}^{n,p,\gamma}(a,h):=\set{t_{x}^{\alpha,p}=t_x^{\alpha,p}\left(s_x^{\alpha \gamma} \right) \ : \ s_x^{\alpha \gamma}\in {}_m\Ops_{x,\alpha}^{n,\gamma}(h) \ \mbox{ and } \ t_{x}^{\alpha,p} h(x):=    \sum_{\abs{\gamma} =0}^p \dfrac{1}{\gamma !}\hat{e}_j s_x^{\alpha \gamma} [h]_j(a)(x-a)^\gamma  +   r_{\alpha\gamma}^p(x-a)  },
\end{eqnarray}

\begin{eqnarray}
{}_m\Opt_{x,\alpha}^{\infty,\gamma}(a,h):=\set{t_{x}^{\alpha,\infty}=t_x^{\alpha,\infty}\left(s_x^{\alpha \gamma} \right) \ : \ s_x^{\alpha \gamma}\in {}_m\Ops_{x,\alpha}^{\infty,\gamma}(h) \ \mbox{ and } \ t_{x}^{\alpha,\infty} h(x):=    \sum_{\abs{\gamma} =0}^\infty \dfrac{1}{\gamma !}\hat{e}_j s_x^{\alpha \gamma} [h]_j(a)(x-a)^\gamma       },
\end{eqnarray}

which allow to generalize the Taylor series expansion of a vector-valued function in multi-index notation \cite{torres2021fracnewrapaitken}. As a consequence, it is possible to obtain the following results:

\begin{align}
 \mbox{If } t_{x}^{\alpha,p}\in {}_m\Opt_{x,\alpha}^{1,p,\gamma}(a,h)  \mbox{ and }  \alpha \to 1  \ \Rightarrow \ &  t_{x}^{1,p} h(x) =  h(a)+  \sum_{\abs{\gamma} =1}^p \dfrac{1}{\gamma !}\hat{e}_j \der{\partial}{x}{\gamma} [h]_j(a)   (x-a)^\gamma   + r_{\gamma}^p(x-a),  \\
\mbox{If } t_{x}^{\alpha,p}\in {}_m\Opt_{x,\alpha}^{1,p,\gamma}(a,h)\mbox{ and } p\to 1 \ \Rightarrow \ &   t_{x}^{\alpha,1} h(x)=   h(a)+  \sum_{k=1}^m \hat{e}_j o_k^{\alpha} [h]_j(a)  \left[ (x-a)\right]_k  + r_{\alpha\gamma }^1(x-a).
\end{align}

Let $f:\Omega \subset \nset{R}^n \to \nset{R}^n$ be a function with a point $\xi \in \Omega$ such that $\norm{f(\xi)}=0$. So, for some $x_i\in B(\xi;\delta)\subset \Omega$ and for some fractional operator $t_{x}^{\alpha,\infty}\in {}_n\Opt_{x,\alpha}^{\infty,\gamma}(x_i,f)$, it is possible to define a type of linear approximation of the function $f$ around the value $x_i$ as follows

\begin{eqnarray*}
t_{x}^{\alpha,\infty} f(x)\approx   f(x_i)+  \sum_{k=1}^n \hat{e}_j o_k^{\alpha} [f]_j(x_i)  \left[ (x-x_i)\right]_k, 
\end{eqnarray*}

which may be rewritten more compactly as follows

\begin{eqnarray}
t_{x}^{\alpha,\infty} f(x) \approx f(x_i)+\left( o_k^\alpha [f]_j(x_i)\right) (x-x_i).
\end{eqnarray}

where $\left( o_k^\alpha [f]_j(x_i)\right)$ denotes a square matrix. On the other hand, if $x\to \xi$ since $\norm{f(\xi)} = 0$ it follows that

\begin{eqnarray*}
0 \approx f(x_i)+\left( o_k^\alpha [f]_j(x_i)\right)(\xi-x_i) & \Rightarrow & \xi \approx x_i-\left( o_k^\alpha [f]_j(x_i)\right)^{-1} f(x_i),
\end{eqnarray*}

then defining the following matrix

\begin{eqnarray}
A_{f,\alpha}(x)=\left( [A_{f,\alpha}]_{jk}(x)\right):=\left(o_k^\alpha[f]_j(x) \right)^{-1},
\end{eqnarray}

it is possible to define the following \textbf{fractional fixed point method}

\begin{eqnarray}
x_{i+1}:=\Phi(\alpha,x_i)=x_i-A_{f,\alpha}(x_i)f(x_i), & i=0,1,2,\cdots,
\end{eqnarray}

which corresponds to the more general case of the \textbf{fractional Newton-Raphson method} \cite{torres2021fracnewrap,torres2021fracnewrapaitken,torres2020fracsome}. As a consequence, considering an iteration function $\Phi:(\nset{R}\setminus \nset{Z})\times \nset{R}^n\to \nset{R}^n$, the iteration function of a fractional iterative method may be written in general form as follows

\begin{eqnarray}\label{eq:4-000}
\Phi(\alpha,x):=x-A_{g,\alpha}(x)f(x),& \alpha\in \nset{R}\setminus \nset{Z},
\end{eqnarray}

where $A_{g,\alpha}$ is a matrix that depends, in at least one of its entries, on fractional operators of order $\alpha$ applied to some function $g:\nset{R}^n \to \nset{R}^n$, whose particular case occurs when $g=f$. So, it is possible to define in a general way a fractional fixed point method as follows

\begin{eqnarray}\label{eq:4-001}
x_{i+1}:=\Phi(\alpha,x_i), & i=0,1,2,\cdots.
\end{eqnarray}

Before continuing, it is necessary to mention that one of the main advantages of fractional iterative methods is that the initial condition $x_0$ can remain fixed, with which it is enough to vary the order $\alpha$ of the fractional operators involved until generating a sequence convergent $\set{x_i}_{i\geq 1}$ to the value $\xi\in \Omega$. Since the order $\alpha$ of the fractional operators is varied, different values of $\alpha$ can generate different sequences convergent to the value $\xi$ but with a different number of iterations (see Figure \ref{fig:4-001}). So, it is possible to define the following set

\begin{eqnarray}
\Conv(\xi):=\set{\Phi \ : \ \lim_{x\to \xi}\Phi(\alpha,x)=x_\xi\in B(\xi;\epsilon)},
\end{eqnarray}

which may be interpreted as the set of fractional fixed point methods that define a convergent sequence $\set{x_i}_{i\geq 1}$ to some value $x_\xi\in B(\xi;\epsilon)$.

\begin{figure}[!ht]
\centering
\includegraphics[width=0.55\textwidth, height=0.35\textwidth]{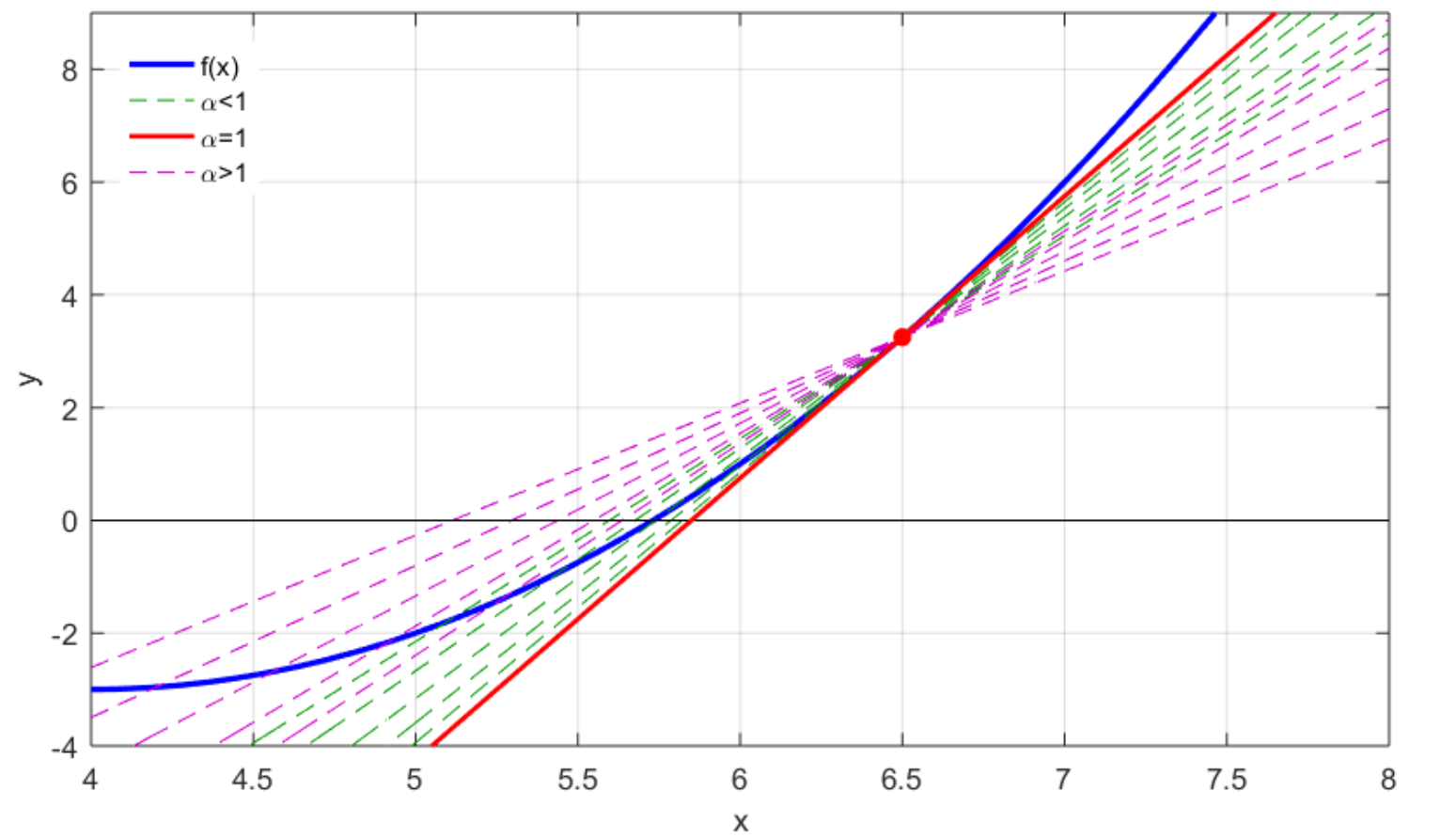}
\caption{Illustration of some lines generated by the fractional Newton-Raphson method for the same initial condition $x_0$ but with different orders $\alpha$ of the fractional operator implemented \cite{torres2021fracnewrap}.}\label{fig:4-001}
\end{figure}

So, denoting by $\card\left(  \cdot  \right)$ the cardinality of a set, it is possible to define the following theorem:

\begin{theorem}\label{teo:4-001}
Let $\Phi:(\nset{R}\setminus \nset{Z})\times\nset{R}^n \to \nset{R}^n$ be an iteration function  with a value $\alpha\in \nset{R}\setminus \nset{Z}$ such that $\Phi(\alpha,x)\in \Conv(\xi)$ in a region $\Omega$. So, if there exists $\epsilon>0$ small enough to warrant that there exists a non-integer value $\beta\in B(\alpha;\epsilon)$ such that

\begin{eqnarray*}
\norm{\Phi(\alpha,x)-\Phi(\beta,x)}<\delta \ \forall x\in \Omega & \mbox{ and } & \Phi(\beta,x)\in \Conv(\xi),
\end{eqnarray*}

then it is fulfilled that

\begin{eqnarray}
\card\left(\Conv(\xi) \right)=\card\left(\nset{R} \right).
\end{eqnarray}

\begin{proof}

The proof of the theorem is carried out by contradiction. Assuming that 

\begin{eqnarray*}
\card\left(\Conv(\xi) \right)<\card\left(\nset{R}\right).
\end{eqnarray*}

So, since $\Phi(\alpha,x),\Phi(\beta,x)\in \Conv(\xi)$ , there exists at least one value $x_k\in B(\xi;\delta)$ such that

\begin{eqnarray*}
\norm{\Phi(\alpha,x_k)-\Phi(\beta,x_k)}=\norm{x_{k+1}-x_{\beta,k+1}}\leq \epsilon_\beta,
\end{eqnarray*}

since $\beta \in B(\alpha ; \epsilon)$ for some $\epsilon$ small enough, without loss of generality if $(n-1)<\alpha < \beta\leq n$ with $n= \lceil \alpha \rceil$,  it follows that

\begin{eqnarray}
\norm{\Phi(\alpha,x_k)-\Phi(a,x_k)}=\norm{x_{k+1}-x_{a,k+1}}\leq \epsilon_a \ \forall a\in[\alpha,\beta],
\end{eqnarray}

as a consequence

\begin{eqnarray*}
\Conv(\xi)\supset \set{\Phi(a,x) \ : \ a \in [\alpha, \beta]} & \Rightarrow & \card \left( \Conv(\xi)\right)> \card \left([\alpha,\beta] \right), 
\end{eqnarray*}

then considering the following function 

\begin{eqnarray*}
h(x)=\dfrac{x-\alpha}{\beta-\alpha},
\end{eqnarray*}

it is fulfilled that

\begin{eqnarray*}
h: [\alpha,\beta]\to[0,1] & \Rightarrow &  \card \left([\alpha,\beta] \right)=\card \left([0,1] \right)=\card\left(\nset{R} \right),
\end{eqnarray*}

and therefore

\begin{eqnarray*}
 \card \left( \Conv(\xi)\right)> \card\left(\nset{R} \right).
\end{eqnarray*}

\end{proof}

\end{theorem}

Finally, it is necessary to mention that fractional iterative methods may be defined in the complex space \cite{torres2021fracnewrapaitken,torres2021zeta}, that is,

\begin{eqnarray}
\set{\Phi(\alpha,x) \ : \ \alpha \in \nset{R}\setminus \nset{Z} \ \mbox{ and } \ x\in \nset{C}^n  }.
\end{eqnarray}

However, due to the part of the integral operator that fractional operators usually have, it may be considered that in the matrix $A_{g,\alpha}$ each fractional operator $o_k^\alpha$ is obtained for a real variable $[x]_k$, and if the result allows it, this variable is subsequently substituted by a complex variable $[x_i]_k$, that is,

\begin{eqnarray}
A_{g,\alpha}(x_i):=A_{g,\alpha}(x)\bigg{|}_{ x\longrightarrow x_i}  , & x\in \nset{R}^n, & x_i\in \nset{C}^n.
\end{eqnarray}

So, considering the above as well as the \mref{Theorem}{\ref{teo:2-001}} and the \mref{Theorem}{\ref{teo:4-001}}, the following corollary is obtained:

\begin{corollary}
Let $\Phi: (\nset{R} \setminus \nset{Z})\times \nset{C}^n \to \nset{C}^n$ be an interaction function with a sequence of different values $\set{\alpha_i}_{i\geq 1} \in \nset{R}\setminus \nset{Z}$ such that defines the following set

\begin{eqnarray*}
\Conv\left(\Omega, \set{\alpha_i}_{i\geq 1} \right):=\set{\Phi(\alpha,x)\in \Conv\left(\xi_\alpha \right) \mbox{ for some } \xi_\alpha\in \Omega \ : \ \alpha \in \set{\alpha_i}_{i\geq 1} }.
\end{eqnarray*}

So, if $\card\left(\Conv\left(\Omega, \set{\alpha_i}_{i\geq 1} \right)\right)=M$ with $1<M<\infty$, then $\Phi$ has a mean order of convergence of order (at least) $\overline{p}$ in $\Omega$, and there exists a sequence  $\set{P_{i,\alpha_i}}_{i\geq 1}^M\in B(\overline{p};\delta_K) $ that allows defining the following value

\begin{eqnarray*}
\overline{P}:=\dfrac{1}{M}\sum_{i=1}^M P_{i,\alpha_i},
\end{eqnarray*}

and therefore, for $M$ large enough it is fulfilled that

\begin{eqnarray}
\abs{\overline{P}-\overline{p}}<\epsilon.
\end{eqnarray}

\end{corollary}

\section{Approximation to the Critical Points of a Function}

Let $C^s(\Omega)$ a set of functions defined as follows

\begin{eqnarray}
C^s(\Omega):= \set{f \ : \ \exists \der{\partial}{x}{\gamma}f(x) \ \forall \abs{\gamma}\leq s \ \mbox{ and } \ \forall x\in \Omega  }.
\end{eqnarray}

So, it is possible to obtain the following result:

\begin{eqnarray}
\mbox{Let }f:\Omega \subset \nset{R}^n \to \nset{R} \mbox{ a function such that }  f\in C^2(\Omega) \ \Rightarrow \ \exists \nabla f(x) \mbox{ and }\exists Hf(x) \ \forall x\in \Omega,
\end{eqnarray}

where $\nabla f$ and $Hf$ denote the gradient of $f$ and the Hessian matrix of $f$ respectively. So in general, for every scalar function $f: \Omega \subset \nset{C}^n \to \nset{C}$ that belongs to the set $C^2(\Omega)$, it is possible to define the following set

\begin{eqnarray}
\mathfrak{C}(\Omega,f):=  \set{\xi \in \Omega \ : \ \norm{\nabla f(\xi)}=0},
\end{eqnarray}

which corresponds to the set of critical points of the function $f$ in the set $\Omega$. On the other hand, denoting by $\re\left(\cdot \right)$ the real part of a complex, by $\det\left(\cdot \right)$ the determinant of a matrix and  by $\sgn \left( \cdot \right)$ the sign function such that for a square matrix $A$

\begin{eqnarray*}
\sgn\left(A \right):= \left( \sgn\left([A]_{jk} \right) \right),
\end{eqnarray*}

it is possible to define the following functions

\begin{eqnarray}
\Delta_{d}(\xi):= \sgn\left(\det\left(\re\left(Hf(\xi) \right) \right)\right)& \mbox{ and }&
\Delta_{t}(\xi):= \tr\left(\sgn\left(\re\left(Hf(\xi) \right)\right) \right),
\end{eqnarray}

which allow to define the following sets

\begin{eqnarray}
\mathfrak{C}_M(\Omega,f):= \set{\xi \in \mathfrak{C}(\Omega,f) \ : \ \Delta_{d}(\xi)=1  \ \mbox{ and } \ \Delta_{t}(\xi)=-n   },
\end{eqnarray}

\begin{eqnarray}
\mathfrak{C}_m(\Omega,f):= \set{\xi \in \mathfrak{C}(\Omega,f) \ : \ \Delta_{d}(\xi)=1  \ \mbox{ and } \ \Delta_{t}(\xi)=n  },
\end{eqnarray}

\begin{eqnarray}
\mathfrak{C}_S(\Omega,f):= \set{\xi \in \mathfrak{C}(\Omega,f) \ : \ \Delta_{d}(\xi)=-1  \ \mbox{ and } \ \Delta_{t}(\xi)\in [-n,n]   },
\end{eqnarray}

which correspond respectively to the sets of local maxima, local minima, and local saddle points of the function $f$ in the set $\Omega$. So, defining the following set of functions

\begin{eqnarray}
C^2_H(\Omega):= \set{f\in C^2(\Omega) \ : \ \exists \left(Hf(x) \right)^{-1}  \ \forall x\in \Omega  },
\end{eqnarray}

and considering a function $f: \Omega \subset \nset{C}^n \to \nset{C}^n $ such that $f\in C^2_H(\Omega)$, it is possible to construct an iteration function $\Phi_{H,\delta}:(\nset{R}\setminus \nset{Z})\times \nset{C}^n \to \nset{C}^n$ defined as follows

\begin{eqnarray}
\Phi_{H,\delta}(\alpha,x):=x-\mathcal{H}_{g,\alpha}(x)   \nabla f(x),
\end{eqnarray}

which corresponds to the iteration function of a hybrid fractional iterative method, where 

\begin{eqnarray}
\mathcal{H}_{g,\alpha}(x):=\left\{
\begin{array}{cc}
A_{g,\alpha}(x), & \mbox{If } \norm{\nabla f(x)}>\delta \vspace{0.1cm}\\
\left(Hf(x) \right)^{-1}, & \mbox{If }\norm{\nabla f(x)}\leq \delta
\end{array}\right. ,
\end{eqnarray}

and $A_{g,\alpha}$  is a matrix of some fractional fixed point method.

\subsection{Examples}

Let $f:\Omega \subset \nset{C}^2 \to \nset{C} $ be a function given by the following expression

\begin{eqnarray}
f(x)=(2 - x^2 + x^3 y) cos(x) - (2 - y^2) cos(y) - x (5 - y^3 cos(y) - 2 sin(x)) - y (7 + 2 sin(y)),
\end{eqnarray}

then

\begin{eqnarray}
\nabla f(x)=\begin{pmatrix}
3 x^2 y cos(x) + y^3 cos(y) + x^2 (1 - x y) sin(x)-5  \vspace{0.1cm}\\ 
x^3 cos(x) + 3 x y^2 cos(y) - y^2(1+ x y)sin(y)-7 
\end{pmatrix},
\end{eqnarray}

\begin{eqnarray}
H f(x)=\begin{pmatrix}
x ((x + 6 y - x^2 y) cos(x) + 2 (1 - 3 x y) sin(x)) & x^2(3 cos(x)- x sin(x)) + y^2(3  cos(y)  - y sin(y)) \vspace{0.1cm}\\ 
x^2(3 cos(x) - x sin(x) )+ y^2(3 cos(y) - ysin(y)) & -y ((y - x (6-y^2)) cos(y) + 2 (1 + 3 x y) sin(y))
\end{pmatrix}.
\end{eqnarray}

So, considering the following function

\begin{eqnarray}\label{eq:5-000}
\rnd_m\left(x \right):=\left\{
\begin{array}{cc}
\re\left(x\right),& \mbox{if } \abs{\im(x)}\leq 10^{-m}\vspace{0.1cm}\\
x,& \mbox{if } \abs{\im(x)}> 10^{-m}\vspace{0.1cm}\\
\end{array}\right., 
\end{eqnarray}

it is possible to define the following iteration function

\begin{eqnarray}\label{eq:5-001}
\rnd_5\left(\Phi_{H,\delta}(\alpha,x) \right):=\hat{e}_j\rnd_5\left(\left[\Phi_{H,\delta} \right]_j(\alpha,x)\right).
\end{eqnarray}

Before continuing, it is necessary to mention that a description of the algorithm that must be implemented when working with a fractional iterative method given by the equation \eqref{eq:4-001} may be found in the reference \cite{torres2020fracsome}. Simplified examples of how a fractional iterative method given by a matrix $A_{g,\alpha}$ should be programmed may be found in the references \cite{torres2021codefracquasi,torres2021codeaccelfracquasi,torres2021codefracpseudo}.

\begin{example}
Using the function \eqref{eq:5-000}, the Riemann-Liouville fractional derivative \eqref{eq:3-002} and $\nabla f$, it is possible to construct an iteration function analogous to the equation \eqref{eq:4-000} using the following matrix

\begin{eqnarray}\label{eq:5-002}
A_{g,\alpha}(x_i)=A_{g_f,\beta}(x_i)=\left([A_{g_f,\beta}]_{jk}(x_i) \right) :=\left( \partial_k^{\beta(\alpha,[x]_k)}[g_f]_{j}(x)  \right)_{x_i}^{-1}, & \alpha\in \nset{R}\setminus \nset{Z},
\end{eqnarray}

which generates a particular case of the \textbf{fractional quasi-Newton method} \cite{torres2020fracsome,torres2021acceleration}, where $ g_{f}(x) $  and $ \beta (\alpha, [x]_k) $ are functions defined as follows

\begin{eqnarray}\label{eq:5-003}
g_{f}(x):=\nabla f(x_i)+Hf(x_i)x& \mbox{ and }
&
\beta(\alpha,[x]_k):=\left\{
\begin{array}{cc}
\alpha, &\mbox{if \hspace{0.1cm} } \abs{  [x]_k }\neq 0 \vspace{0.1cm}\\
1,& \mbox{if \hspace{0.1cm}  }  \abs{  [x]_k  }=0
\end{array}\right. .
\end{eqnarray}

So, considering following initial condition,

\begin{eqnarray*}
x_0=(5.21,5.21)^T  &\mbox{ with } &\norm{\nabla f(x_0)} \approx 1,289.4083,
\end{eqnarray*}

the following results are obtained:

\begin{scriptsize}
\centering
\begin{longtable}{c|ccccccccc}
\toprule
&$\alpha$& $[ x_n ]_1$ &$[ x_n ]_2$ &$\norm{x_n-x_{n-1}}_2$&$\norm{\nabla f(x_n)}_2$&$P_n$&$\Delta_{d}(x_n)$&$\Delta_{t}(x_n)$ &$n$\\
\midrule
    1     & -0.530515 & 6.6771554 - 0.02130862i &  -0.014023 + 1.72836829i & 1.41E-08 & 9.24E-05 & 0.9812 & -1    & 2     & 167 \\
    2     & -0.516037 & 0.01499973 - 1.73190718i &  6.6757499 - 0.04157569i & 1.41E-08 & 9.86E-05 & 1.0260 & -1    & -2    & 165 \\
    3     & -0.472867 & 0.01499966 + 1.73190711i &  6.67574974 + 0.04157578i & 2.45E-08 & 9.54E-05 & 1.0000 & -1    & -2    & 180 \\
    4     & -0.440017 & 6.67715551 + 0.02130861i &  -0.014023 - 1.72836833i & 1.41E-08 & 9.47E-05 & 1.0113 & -1    & 2     & 180 \\
    5     & -0.372536 & -1.12922862 + 1.02480512i &  3.7817693 + 0.02894647i & 3.22E-07 & 8.23E-05 & 0.9960 & 1     & -2    & 92 \\
    6     & -0.359168 & -1.12922793 - 1.02480539i &  3.78176969 - 0.02894643i & 1.36E-06 & 7.78E-05 & 1.0255 & 1     & -2    & 92 \\
    7     & -0.317767 & 3.68514423 - 0.05398726i &  -1.20114465 + 1.03004598i & 4.35E-07 & 7.98E-05 & 1.0095 & 1     & -2    & 89 \\
    8     & -0.175657 & 6.66385192 + 0.00958153i &  -3.05535188 + 0.51526774i & 1.66E-07 & 9.98E-05 & 1.0145 & 1     & 2     & 129 \\
    9     & -0.174937 & 9.69844564 - 0.00485976i &  -1.49692201 + 1.85490018i & 1.41E-08 & 8.99E-05 & 0.9812 & 1     & -2    & 180 \\
    10    & -0.167409 & 9.69844566 + 0.00485981i &  -1.49692201 - 1.85490019i & 1.00E-08 & 8.76E-05 & 1.0000 & 1     & -2    & 178 \\
    11    & -0.165538 & 3.68514454 + 0.0539876i &  -1.20114479 - 1.0300467i & 8.66E-07 & 8.57E-05 & 1.0144 & 1     & -2    & 117 \\
    12    & -0.162111 & -1.47430587 + 1.85378122i &  9.71215809 + 0.012692i & 1.41E-08 & 8.26E-05 & 1.0313 & 1     & -2    & 178 \\
    13    & -0.148486 & 12.78190313 - 0.00664448i &  -3.36083258 - 1.47015693i & 1.00E-08 & 8.71E-05 & 1.0192 & 1     & 2     & 195 \\
    14    & -0.141354 & -1.47430585 - 1.85378123i &  9.71215813 - 0.01269197i & 3.00E-08 & 5.73E-05 & 0.9966 & 1     & -2    & 179 \\
    15    & -0.140788 & -3.01831349 + 0.5058919i &  6.69924174 + 0.01613682i & 3.16E-08 & 9.57E-05 & 1.0285 & 1     & 2     & 146 \\
    16    & -0.125015 & 19.0075656 & -7.54961078 & 1.41E-08 & 8.38E-05 & 1.0000 & 1     & 2     & 197 \\
    17    & -0.119655 & -4.59285859 & 9.73129666 & 3.61E-08 & 4.16E-05 & 1.0195 & 1     & -2    & 111 \\
    18    & -0.092015 & 6.66385199 - 0.00958166i &  -3.05535203 - 0.51526774i & 2.45E-08 & 8.85E-05 & 1.0044 & 1     & 2     & 85 \\
    19    & -0.081244 & 12.81002482 & -7.10966547 & 1.41E-08 & 8.10E-05 & 1.0399 & 1     & 2     & 190 \\
    20    & -0.075076 & 9.71878342 & -4.62771758 & 2.24E-08 & 9.26E-05 & 1.0000 & 1     & -2    & 97 \\
    21    & -0.073120 & -3.01831348 - 0.50589194i &  6.69924174 - 0.01613673i & 4.58E-08 & 7.95E-05 & 1.0187 & 1     & 2     & 82 \\
    22    & -0.056190 & 18.99311678 & -9.30049381 & 1.41E-08 & 9.76E-05 & 0.9812 & -1    & 0     & 145 \\
    23    & -0.052492 & -6.39937485 & 9.68519629 & 2.24E-08 & 9.90E-05 & 0.9563 & -1    & 0     & 161 \\
    24    & -0.052490 & -7.09665187 & 12.81542466 & 2.24E-08 & 7.76E-05 & 1.0377 & 1     & 2     & 113 \\
    25    & -0.037197 & -5.68870793 - 0.65962195i &  15.8889979 - 0.00516137i & 1.41E-08 & 2.87E-05 & 0.9812 & 1     & -2    & 183 \\
    26    & -0.030387 & -9.30202535 & 18.99474019 & 1.41E-08 & 7.22E-05 & 0.9812 & -1    & 0     & 162 \\
 \bottomrule
\caption{Results obtained using the fractional quasi-Newton method \cite{torres2021acceleration}.}
\end{longtable}
\end{scriptsize}

Therefore 

\begin{eqnarray}
\overline{P}\approx 1.0060\in B\left(\overline{p};\delta_K \right), 
\end{eqnarray}

which is consistent with the \mref{Corollary}{\ref{cor:2-001}}, since in general if $\xi \in \mathfrak{C}(\Omega,f)$, then it is fulfilled that \cite{torres2021fracnewrapaitken}

\begin{eqnarray*}
\lim_{x\to \xi}\norm{\Phi^{(1)}(\alpha,x)}\neq 0.
\end{eqnarray*}

\end{example}

\begin{example}
Using the iteration function \eqref{eq:5-001} and the matrix $A_{g,\alpha}$ given by the equation  \eqref{eq:5-002}, considering following values

\begin{eqnarray*}
\delta=7 & \mbox{ and } &   x_0=(4.78,4.78)^T  \ \mbox{ with } \  \norm{\nabla f(x_0)} \approx 770.4734,
\end{eqnarray*}

the following results are obtained:

\begin{scriptsize}
\centering
\begin{longtable}{c|ccccccccc}
\toprule
&$\alpha$& $[ x_n ]_1$ &$[ x_n ]_2$ &$\norm{x_n-x_{n-1}}_2$&$\norm{\nabla f(x_n)}_2$&$P_n$&$\Delta_{d}(x_n)$&$\Delta_{t}(x_n)$ &$n$\\
\midrule
    1     & -0.991504 & 3.98115471 & 3.92170125 & 1.00E-08 & 1.50E-06 & 2.1162 & 1     & 2     & 55 \\
    2     & -0.985320 & -0.20172521 & -2.13862013 & 1.00E-08 & 3.55E-08 & 2.0096 & -1    & 0     & 184 \\
    3     & -0.977534 & 4.76944744 & 0.24682585 & 5.21E-06 & 4.75E-07 & 2.0536 & -1    & 2     & 115 \\
    4     & -0.957378 & -0.14249533 & 7.84459109 & 2.32E-06 & 1.71E-06 & 2.1574 & -1    & 0     & 44 \\
    5     & -0.931674 &  1.52183063 + 0.04852431i &  -1.07285283 + 0.62177498i & 1.64E-05 & 6.40E-08 & 1.9728 & -1    & 0     & 147 \\
    6     & -0.910766 & -0.1411895 & 4.75629836 & 5.06E-07 & 4.42E-07 & 2.0939 & -1    & -2    & 99 \\
    7     & -0.902424 & -1.66983169 & -1.47843397 & 3.51E-05 & 6.06E-08 & 2.0210 & 1     & -2    & 141 \\
    8     & -0.796926 & 7.84012182 & 0.11780088 & 5.96E-05 & 2.18E-06 & 2.2020 & -1    & 0     & 32 \\
    9     & -0.747172 &  -1.47430586 - 1.85378123i &  9.71215811 - 0.01269197i & 5.21E-06 & 1.45E-05 & 2.0987 & 1     & -2    & 193 \\
    10    & -0.739854 &  9.69844563 - 0.0048598i &  -1.496922 + 1.85490017i & 4.57E-06 & 1.59E-05 & 2.1076 & 1     & -2    & 190 \\
    11    & -0.734400 &  9.69844563 + 0.0048598i &  -1.496922 - 1.85490017i & 4.77E-06 & 1.59E-05 & 2.1051 & 1     & -2    & 194 \\
    12    & -0.718024 &  -1.47430586 + 1.85378123i &  9.71215811 + 0.01269197i & 5.09E-06 & 1.45E-05 & 2.1055 & 1     & -2    & 172 \\
    13    & -0.691512 &  -1.12922847 - 1.02480556i &  3.78176946 - 0.02894603i & 4.21E-06 & 5.54E-07 & 2.0281 & 1     & -2    & 166 \\
    14    & -0.654774 &  -0.9615658 + 0.5828065i &  1.85727226 + 0.22306481i & 2.46E-06 & 1.41E-07 & 1.9957 & -1    & 0     & 99 \\
    15    & -0.639046 &  0.72967089 + 0.94166299i &  0.62407461 - 0.91988663i & 7.18E-07 & 1.54E-07 & 2.0484 & 1     & 2     & 128 \\
    16    & -0.616404 &  3.68514466 + 0.05398708i &  -1.20114498 - 1.03004629i & 6.54E-06 & 6.96E-07 & 1.9881 & 1     & -2    & 150 \\
    17    & -0.598098 &  -1.12922847 + 1.02480556i &  3.78176946 + 0.02894603i & 3.10E-06 & 5.54E-07 & 2.0471 & 1     & -2    & 62 \\
    18    & -0.591784 &  3.68514466 - 0.05398708i &  -1.20114498 + 1.03004629i & 8.24E-06 & 6.96E-07 & 2.0008 & 1     & -2    & 67 \\
    19    & -0.531176 &  6.67715546 - 0.02130875i &  -0.01402295 + 1.7283683i & 1.41E-08 & 4.70E-06 & 1.9773 & -1    & 2     & 52 \\
    20    & -0.527738 &  12.78275364 - 0.00603578i &  -0.00730626 + 2.36240058i & 8.25E-07 & 2.69E-05 & 2.1144 & -1    & 2     & 193 \\
    21    & -0.511182 &  1.59511265 + 0.92462709i &  0.28169602 - 0.00845802i & 3.61E-08 & 9.30E-08 & 1.9993 & -1    & 2     & 70 \\
    22    & -0.503186 &  0.01499973 - 1.73190712i &  6.67574976 - 0.04157565i & 3.33E-05 & 3.96E-06 & 2.1931 & -1    & -2    & 57 \\
    23    & -0.490941 &  -3.34309333 + 1.46646036i &  12.79048871 + 0.01073275i & 7.06E-07 & 4.43E-05 & 2.1248 & 1     & 2     & 194 \\
    24    & -0.490753 &  0.00737884 - 2.36289538i &  12.78266688 - 0.00836806i & 8.19E-07 & 3.00E-05 & 2.1172 & -1    & -2    & 199 \\
    25    & -0.470183 &  12.78275364 + 0.00603578i &  -0.00730626 - 2.36240058i & 8.35E-07 & 2.69E-05 & 2.1169 & -1    & 2     & 200 \\
    26    & -0.468001 &  -3.34309333 - 1.46646036i &  12.79048871 - 0.01073275i & 9.49E-07 & 4.43E-05 & 2.0622 & 1     & 2     & 186 \\
    27    & -0.463959 &  12.78190312 - 0.00664448i &  -3.36083257 - 1.47015693i & 3.42E-07 & 3.36E-05 & 2.1539 & 1     & 2     & 200 \\
    28    & -0.458777 &  1.30993837 - 0.36023537i &  0.99738945 - 0.66890573i & 1.16E-06 & 8.98E-08 & 1.9828 & -1    & 2     & 53 \\
    29    & -0.437585 &  0.01499973 + 1.73190712i &  6.67574975 + 0.04157565i & 8.14E-05 & 5.27E-06 & 2.0388 & -1    & -2    & 57 \\
    30    & -0.429119 &  12.78190312 + 0.00664448i &  -3.36083257 + 1.47015693i & 2.80E-07 & 3.36E-05 & 2.1486 & 1     & 2     & 184 \\
    31    & -0.417531 &  6.67715546 + 0.02130875i &  -0.01402295 - 1.72836831i & 8.14E-05 & 5.37E-06 & 2.0768 & -1    & 2     & 49 \\
    32    & -0.321303 &  15.88192661 + 0.00033296i &  -1.64153442 - 2.37001819i & 1.37E-07 & 4.16E-05 & 2.1186 & 1     & -2    & 192 \\
    33    & -0.295259 &  15.88518055 + 0.00474592i &  -5.70013516 + 0.67487422i & 4.69E-08 & 5.96E-05 & 2.1502 & 1     & -2    & 195 \\
    34    & -0.287905 &  -5.68870793 + 0.65962195i &  15.8889979 + 0.00516137i & 7.23E-07 & 2.87E-05 & 2.0120 & 1     & -2    & 177 \\
    35    & -0.278601 &  15.88518055 - 0.00474592i &  -5.70013516 - 0.67487422i & 1.15E-07 & 5.96E-05 & 2.0524 & 1     & -2    & 197 \\
    36    & -0.264047 &  -5.68870793 - 0.65962195i &  15.8889979 - 0.00516137i & 3.32E-08 & 2.87E-05 & 2.1731 & 1     & -2    & 194 \\
    37    & -0.263797 &  6.66385192 - 0.00958162i &  -3.05535199 - 0.51526776i & 3.78E-05 & 5.76E-06 & 2.0466 & 1     & 2     & 110 \\
    38    & -0.242447 & -4.59285856 & 9.73129667 & 1.34E-07 & 2.16E-05 & 2.7985 & 1     & -2    & 199 \\
    39    & -0.240107 & 9.71878344 & -4.6277176 & 5.48E-07 & 7.22E-06 & 2.6624 & 1     & -2    & 173 \\
    40    & -0.235095 &  -3.01831353 + 0.50589193i &  6.69924181 + 0.01613676i & 3.61E-08 & 1.43E-06 & 1.9762 & 1     & 2     & 77 \\
    41    & -0.212867 &  6.66385192 + 0.00958162i &  -3.05535199 + 0.51526775i & 8.78E-05 & 6.78E-06 & 1.9815 & 1     & 2     & 57 \\
    42    & -0.211725 & 19.0075656 & -7.54961079 & 1.61E-04 & 4.83E-05 & 0.7767 & 1     & 2     & 197 \\
    43    & -0.209337 &  -3.01831353 - 0.50589194i &  6.69924181 - 0.01613676i & 7.73E-05 & 3.05E-06 & 1.9919 & 1     & 2     & 64 \\
    44    & -0.204931 & -7.53686364 & 19.00985885 & 1.00E-08 & 3.66E-05 & 2.0867 & 1     & 2     & 158 \\
    45    & -0.181783 & 12.81002482 & -7.10966546 & 1.32E-07 & 3.93E-05 & 2.2517 & 1     & 2     & 196 \\
    46    & -0.181407 & -9.30202535 & 18.99474019 & 1.00E-08 & 7.22E-05 & 2.1044 & -1    & 0     & 197 \\
    47    & -0.178655 & -7.09665188 & 12.81542466 & 1.00E-08 & 4.79E-05 & 2.6959 & 1     & 2     & 188 \\
    48    & -0.175623 & 18.99311678 & -9.3004938 & 1.00E-08 & 6.54E-05 & 2.1290 & -1    & 0     & 187 \\
    49    & -0.125919 & -6.39937487 & 9.6851963 & 1.81E-06 & 2.11E-05 & 2.0664 & -1    & 0     & 195 \\
    50    & -0.092457 & 9.67778512 & -6.40235748 & 5.02E-07 & 2.22E-05 & 2.2493 & -1    & 0     & 183 \\
    51    & -0.076797 & 19.02754978 & -12.95559618 & 1.29E-04 & 4.95E-05 & 0.9503 & 1     & 2     & 156 \\
 \bottomrule
\caption{Results obtained using the iteration function \eqref{eq:5-001} with the fractional quasi-Newton method \cite{torres2021acceleration}.}
\end{longtable}
\end{scriptsize}

Therefore 

\begin{eqnarray}
\overline{P}\approx 2.0692 \in B\left(\overline{p};\delta_K \right), 
\end{eqnarray}

which is consistent with the \mref{Corollary}{\ref{cor:2-001}}, since in general if $\xi \in \mathfrak{C}(\Omega,f)$, then it is fulfilled that \cite{torres2021fracnewrapaitken}

\begin{eqnarray*}
\lim_{x\to \xi}\norm{\Phi_{H,\delta}^{(1)}(\alpha,x)}=0.
\end{eqnarray*}

\end{example}

\begin{example}

Using the Riemann-Liouville fractional derivative \eqref{eq:3-002}, it is possible to construct the following matrix

\begin{eqnarray}\label{eq:5-004}
A_{g,\alpha}(x_i)=A_{\epsilon,\beta}(x_i)=\left([A_{\epsilon,\beta}]_{jk}(x_i)\right):=\left( \partial_k^{\beta(\alpha,[x]_k)}\delta_{jk}+ \epsilon\delta_{jk}  \right)_{x_i}, & \alpha\in \nset{R}\setminus \nset{Z},
\end{eqnarray}

which generates a particular case of the \textbf{fractional pseudo-Newton method} \cite{torres2020fracpseunew,torres2021zeta,torres2020reduction},  where $ \delta_{jk} $ is the Kronecker delta, $ \epsilon $ is a positive constant $ \ll 1 $, and $ \beta (\alpha, [x_i]_k) $  is a function defined by the equation \eqref{eq:5-003}. So, using the iteration function \eqref{eq:5-001} and the matrix $A_{g,\alpha}$ given by the equation  \eqref{eq:5-004}, considering following values

\begin{eqnarray*}
\begin{array}{cccc}
\epsilon=10^{-4}, &\delta=13 & \mbox{ and } &   x_0=(14.55,14.55)^T \  \mbox{ with } \  \norm{\nabla f(x_0)} \approx 65,057.2221,
\end{array}
\end{eqnarray*}

the following results are obtained:

\begin{scriptsize}
\centering
\begin{longtable}{c|ccccccccc}
\toprule
&$\alpha$& $[ x_n ]_1$ &$[ x_n ]_2$ &$\norm{x_n-x_{n-1}}_2$&$\norm{\nabla f(x_n)}_2$&$P_n$&$\Delta_{d}(x_n)$&$\Delta_{t}(x_n)$ &$n$\\
\midrule
    1     & 0.997025 & 6.40346174 & -9.68745629 & 7.48E-06 & 2.99E-05 & 2.0712 & -1    & 0     & 11 \\
    2     & 0.997053 & -6.8254374 & -6.80736533 & 1.34E-06 & 1.17E-05 & 2.1970 & 1     & -2    & 37 \\
    3     & 0.997061 & 9.73394944 & 4.59418309 & 3.34E-06 & 1.33E-05 & 2.1262 & 1     & 2     & 33 \\
    4     & 0.998113 & 4.62598971 & 9.72138809 & 2.20E-07 & 1.78E-05 & 2.8079 & 1     & 2     & 13 \\
    5     & 0.998133 & -9.67933962 & 6.40821255 & 3.58E-07 & 3.16E-05 & 2.1427 & -1    & 0     & 19 \\
    6     & 0.998185 &  -3.75670368 + 0.00677324i &  1.14479461 - 0.90835133i & 6.32E-08 & 8.42E-07 & 1.9860 & 1     & -2    & 184 \\
    7     & 0.998189 &  -3.75670368 - 0.00677324i &  1.14479461 + 0.90835133i & 2.47E-06 & 8.42E-07 & 1.9809 & 1     & -2    & 126 \\
    8     & 0.998229 & -12.81526848 & -7.09878784 & 3.61E-08 & 3.00E-05 & 2.1703 & 1     & -2    & 22 \\
    9     & 0.998469 & -12.6804252 & -15.85472455 & 1.00E-08 & 4.54E-05 & 2.2093 & -1    & 0     & 49 \\
    10    & 0.999045 &  1.52183063 - 0.04852431i &  -1.07285283 - 0.62177498i & 8.25E-06 & 6.40E-08 & 1.9673 & -1    & 0     & 161 \\
    11    & 0.999065 & 7.09845974 & -12.81449874 & 7.07E-08 & 2.76E-05 & 2.2122 & 1     & 2     & 33 \\
    12    & 0.999909 & 9.81602358 & 9.80895121 & 2.24E-08 & 4.07E-05 & 2.1124 & 1     & 2     & 25 \\
    13    & 0.999917 & -7.09665188 & 12.81542466 & 1.06E-07 & 4.79E-05 & 2.1726 & 1     & 2     & 26 \\
    14    & 0.999921 & -6.80274842 & 6.8263687 & 7.69E-06 & 1.15E-05 & 2.2126 & 1     & 2     & 28 \\
    15    & 0.999925 & -12.80936242 & 7.11220453 & 1.30E-07 & 5.55E-05 & 2.2710 & 1     & 2     & 28 \\
    16    & 0.999929 & -9.73194065 & -4.58368411 & 2.72E-06 & 5.55E-06 & 2.7275 & 1     & 2     & 62 \\
    17    & 0.999937 & -9.81505776 & -9.80760476 & 1.13E-04 & 2.80E-05 & 1.4237 & 1     & 2     & 18 \\
    18    & 0.999941 & -4.61844557 & -9.71852806 & 1.53E-06 & 1.39E-05 & 2.8405 & 1     & 2     & 61 \\
    19    & 0.999945 & 6.80674644 & -6.820744 & 4.50E-06 & 1.44E-05 & 2.1855 & 1     & 2     & 176 \\
    20    & 0.999953 & 12.81002482 & -7.10966546 & 3.41E-07 & 3.93E-05 & 2.2868 & 1     & 2     & 44 \\
    21    & 1.003393 & 6.82167482 & 6.80212518 & 8.31E-06 & 1.48E-05 & 2.1795 & 1     & -2    & 5 \\
    22    & 1.004893 &  -0.55742729 - 0.65679566i &  -0.20882106 - 1.14800938i & 3.11E-07 & 1.41E-07 & 2.0709 & -1    & 0     & 64 \\
    23    & 1.004925 &  3.68514466 + 0.05398708i &  -1.20114498 - 1.03004629i & 5.10E-08 & 6.96E-07 & 1.9686 & 1     & -2    & 119 \\
    24    & 1.004969 &  -1.12922847 + 1.02480556i &  3.78176946 + 0.02894603i & 4.58E-08 & 5.54E-07 & 1.9971 & 1     & -2    & 137 \\
    25    & 1.005025 &  0.72967089 - 0.94166299i &  0.62407461 + 0.91988663i & 2.37E-05 & 1.54E-07 & 1.9983 & 1     & 2     & 84 \\
    26    & 1.005549 & 0.29601303 & -4.65165906 & 1.49E-05 & 4.30E-07 & 2.1087 & -1    & -2    & 15 \\
    27    & 1.005849 &  3.68514466 - 0.05398708i &  -1.20114498 + 1.03004629i & 6.25E-08 & 6.96E-07 & 1.9890 & 1     & -2    & 184 \\
    28    & 1.005937 &  -1.12922847 - 1.02480556i &  3.78176946 - 0.02894603i & 1.41E-08 & 5.54E-07 & 1.9735 & 1     & -2    & 82 \\
    29    & 1.006421 &  -1.3914151 - 0.70003547i &  0.17621271 + 1.00035774i & 1.02E-04 & 1.46E-07 & 2.0270 & -1    & 0     & 50 \\
    30    & 1.006437 &  1.30993837 - 0.36023537i &  0.99738945 - 0.66890573i & 4.91E-06 & 8.98E-08 & 1.9863 & -1    & 2     & 44 \\
    31    & 1.006465 &  -0.55742729 + 0.65679566i &  -0.20882106 + 1.14800938i & 6.32E-08 & 1.41E-07 & 2.1428 & -1    & 0     & 38 \\
    32    & 1.007481 & -3.95538299 & -3.88543329 & 9.14E-05 & 3.64E-06 & 2.3031 & 1     & 2     & 5 \\
    33    & 1.008713 &  1.59511265 - 0.92462709i &  0.28169602 + 0.00845802i & 1.63E-06 & 9.30E-08 & 2.1184 & -1    & 2     & 20 \\
    34    & 1.009697 & -2.30034423 & -0.45950443 & 4.99E-06 & 7.08E-08 & 2.1235 & -1    & 0     & 6 \\
    35    & 1.009817 &  0.09238517 + 0.91135195i &  -1.48626899 - 0.45588717i & 5.70E-07 & 1.37E-07 & 1.9727 & 1     & -2    & 28 \\
    36    & 1.009821 &  0.09238517 - 0.91135195i &  -1.48626899 + 0.45588717i & 1.41E-08 & 1.37E-07 & 2.0053 & 1     & -2    & 34 \\
    37    & 1.009861 &  -1.3914151 + 0.70003546i &  0.17621271 - 1.00035774i & 9.45E-05 & 2.55E-07 & 2.0119 & -1    & 0     & 22 \\
    38    & 1.010385 &  1.30993837 + 0.36023537i &  0.99738945 + 0.66890573i & 4.13E-05 & 8.98E-08 & 1.9803 & -1    & 2     & 38 \\
    39    & 1.908362 &  0.72967089 + 0.94166298i &  0.62407461 - 0.91988663i & 8.45E-05 & 1.83E-07 & 1.9642 & 1     & 2     & 14 \\
    40    & 1.913438 &  1.52183063 + 0.04852431i &  -1.07285283 + 0.62177498i & 1.10E-07 & 6.40E-08 & 1.9787 & -1    & 0     & 13 \\
    41    & 1.918790 & -1.66983169 & -1.47843397 & 1.14E-04 & 6.06E-08 & 2.2493 & 1     & -2    & 5 \\
    42    & 1.920778 &  1.59511265 + 0.92462709i &  0.28169602 - 0.00845802i & 4.58E-08 & 9.30E-08 & 1.9835 & -1    & 2     & 17 \\
    43    & 1.922506 & 3.8890101 & -3.98878888 & 1.22E-07 & 1.48E-06 & 2.1461 & 1     & -2    & 19 \\
    44    & 1.928090 & -3.91843903 & 3.94777085 & 1.97E-07 & 2.03E-06 & 2.0974 & 1     & -2    & 75 \\
    45    & 1.928198 & 4.76944744 & 0.24682585 & 4.45E-05 & 4.75E-07 & 2.0605 & -1    & 2     & 19 \\
    46    & 1.938338 & -0.1411895 & 4.75629836 & 6.65E-06 & 4.42E-07 & 2.0695 & -1    & -2    & 12 \\
    47    & 2.027490 & -4.63811516 & -0.17366027 & 3.50E-07 & 5.12E-07 & 2.4473 & -1    & 2     & 6 \\
    48    & 2.027714 &  -0.9615658 - 0.5828065i &  1.85727226 - 0.22306481i & 1.22E-06 & 1.41E-07 & 2.0016 & -1    & 0     & 80 \\
    49    & 2.027802 &  -0.9615658 + 0.5828065i &  1.85727226 + 0.22306481i & 8.66E-07 & 1.41E-07 & 2.0016 & -1    & 0     & 23 \\
    50    & 2.028082 & 3.98115471 & 3.92170125 & 4.47E-08 & 1.50E-06 & 2.0806 & 1     & 2     & 9 \\
    51    & 2.050222 &  0.10127937 - 0.65790456i &  -0.69552033 - 1.28219351i & 4.24E-08 & 2.55E-08 & 1.9278 & 1     & -2    & 9 \\
    52    & 2.892915 & -0.2017252 & -2.13862013 & 8.96E-05 & 1.79E-07 & 2.0069 & -1    & 0     & 5 \\
    53    & 2.979539 & -9.68548222 & -6.40422387 & 1.48E-05 & 6.43E-06 & 2.0748 & -1    & 0     & 43 \\
    54    & 2.979543 & -6.40734755 & -9.67742959 & 1.68E-05 & 7.33E-06 & 2.0878 & -1    & 0     & 43 \\
    55    & 2.983015 &  6.66385192 - 0.00958162i &  -3.05535199 - 0.51526776i & 4.88E-06 & 5.76E-06 & 1.9701 & 1     & 2     & 65 \\
    56    & 2.983279 &  1.07448447 + 0.94219835i &  -3.88986554 + 0.11532861i & 7.72E-05 & 9.22E-07 & 2.0229 & 1     & -2    & 92 \\
    57    & 2.989991 &  -3.01831353 + 0.50589193i &  6.69924181 + 0.01613676i & 5.64E-06 & 1.43E-06 & 1.9676 & 1     & 2     & 101 \\
    58    & 2.990235 &  12.78190312 + 0.00664448i &  -3.36083257 + 1.47015693i & 3.33E-06 & 3.36E-05 & 2.0443 & 1     & 2     & 27 \\
    59    & 2.990955 &  -3.34309333 - 1.46646036i &  12.79048871 - 0.01073275i & 2.29E-06 & 4.43E-05 & 2.0444 & 1     & 2     & 26 \\
    60    & 3.002283 &  -12.78071432 + 0.00620911i &  3.36250229 + 1.47445201i & 7.91E-06 & 1.73E-05 & 2.0486 & 1     & 2     & 38 \\
    61    & 3.004719 & 9.55477471 & 12.75308268 & 3.30E-07 & 1.49E-05 & 2.1260 & -1    & 0     & 9 \\
    62    & 3.013455 &  -6.65415389 + 0.00918318i &  3.06649242 + 0.56418379i & 5.49E-06 & 4.40E-06 & 1.9795 & 1     & 2     & 90 \\
    63    & 3.013911 & 9.68717241 & 6.39860852 & 1.86E-05 & 7.26E-06 & 2.0743 & -1    & 0     & 199 \\
    64    & 3.014343 & 6.40322967 & 9.6796959 & 2.04E-05 & 1.09E-05 & 2.0718 & -1    & 0     & 189 \\
    65    & 3.982916 &  6.66385192 + 0.00958162i &  -3.05535199 + 0.51526776i & 8.57E-05 & 5.76E-06 & 2.1002 & 1     & 2     & 87 \\
    66    & 3.982992 &  12.78190312 - 0.00664448i &  -3.36083257 - 1.47015693i & 1.25E-05 & 3.36E-05 & 2.0505 & 1     & 2     & 35 \\
    67    & 3.983884 &  3.02691487 + 0.54276524i &  -6.68492207 + 0.01504716i & 1.41E-08 & 5.46E-06 & 1.9751 & 1     & 2     & 117 \\
    68    & 3.990568 &  3.34433054 + 1.46955548i &  -12.78880218 + 0.01004339i & 7.60E-07 & 3.97E-05 & 2.1391 & 1     & 2     & 20 \\
    69    & 3.990580 &  3.34433054 - 1.46955548i &  -12.78880218 - 0.01004339i & 9.72E-07 & 3.97E-05 & 2.1398 & 1     & 2     & 23 \\
    70    & 3.991060 &  -3.01831353 - 0.50589193i &  6.69924181 - 0.01613676i & 9.65E-06 & 1.43E-06 & 1.9672 & 1     & 2     & 81 \\
 \bottomrule
\caption{Results obtained using the iteration function \eqref{eq:5-001} with the fractional psuedo-Newton method \cite{torres2021zeta}.}
\end{longtable}
\end{scriptsize}

Therefore 

\begin{eqnarray}
\overline{P}\approx 2.0994\in B\left(\overline{p};\delta_K \right), 
\end{eqnarray}

which is consistent with the \mref{Corollary}{\ref{cor:2-001}}, since in general if $\xi \in \mathfrak{C}(\Omega,f)$, then it is fulfilled that \cite{torres2021fracnewrapaitken}

\begin{eqnarray*}
\lim_{x\to \xi}\norm{\Phi_{H,\delta}^{(1)}(\alpha,x)}=0.
\end{eqnarray*}

\end{example}

\section{Conclusions}

Considering the large number of fractional operators that exist, and since it does not seem that their number will stop increasing soon at the time of writing this paper, the most simplified and compact way to work the fractional calculus is through the classification of fractional operators using sets which, as shown in the previous sections, allows generalize objects of the conventional calculus such as tensor operators, the diffusion equation, the heat equation, the Taylor series of a vector-valued function and the fixed point method in several variables which allows to generate the method known as the fractional fixed point method that in turn allows a new type of numerical analysis using sets \cite{torres2021acceleration}. Since that each fractional fixed point method that generates a convergent sequence has the ability to generate an uncountable family of fractional fixed point methods that generate convergent sequences as shown by the \mref{Theorem}{\ref{teo:4-001}}, and considering that determining the critical points of a scalar function is usually one of the most recurrent problems in physics, mathematics and engineering, it becomes almost natural to estimate numerically the mean order of convergence of any fractional fixed point method in a region $\Omega$ by determining the critical points of a scalar function. Finally, it should be mentioned that the result of the \mref{Theorem}{\ref{teo:4-001}} can be transferred to the theory of fractional differential equations, resulting in a new type of theory of differential equations using sets, which opens the possibility that the fractional calculus became a more extensive theory and that it should be renamed as \textbf{fractional calculus of sets}.

\bibliography{Biblio}
\bibliographystyle{unsrt}

\end{document}